\documentclass[a4paper,11pt]{amsart}

\usepackage{acatta}

\newcommand{\blankSpace}{\cdot}
\newcommand{\mubar}{\overline{\mu}}

\renewcommand{\Im}{\operatorname{Im}}
\renewcommand{\Re}{\operatorname{Re}}

\title[Almost cplx mflds from the pt of view of Kod dim]{Almost complex manifolds from the point of view of Kodaira dimension}

\author{Andrea Cattaneo}
\address{Dipartimento di Scienze Matematiche, Fisiche e Informatiche
Unit\`a di Matematica e Informatica\\
Universit\`a degli Studi di Parma\\
Parco Area delle Scienze 53/A, 43124\\
Parma, Italy}
\email{andrea.cattaneo@unipr.it}

\thanks{The author is a member of GNSAGA of INdAM. This work is based on a talk given by the author on July 6$^{\text{th}}$ 2021 during the workshop `Cohomology of complex manifolds and special structures II'. He wants to gratefully acknowledge the organizers of the workshop for the kind invitation. The author also thanks the anonymous referee for their review.}

\begin{document}

\begin{abstract}
In complex geometry a classical and useful invariant of a complex manifold is its Kodaira dimension. Since its introduction by Iitaka in the early 70's, its behavior under deformations was object of study and it is known that Kodaira dimension is invariant under holomorphic deformations for a smooth projective manifolds, while there are examples of holomorphic deformations of non-projective manifolds for which the Kodaira dimension is non-constant. Recently this concept has been generalized to almost complex manifolds, we want to present here some of its main features in the non-integrable case, mainly with respect to deformations. At the end we conclude with some speculations on the theory of meromorphic functions on almost complex manifolds.
\end{abstract}

\maketitle

\section{Introduction}

The subject of Kodaira dimension is quite classical in algebraic and complex geometry, which reflects how fast the dimension of the spaces of holomorphic pluricanonical sections grow. A first, rough, classification of complex manifolds can be given in terms of Kodaira dimension. Let us recall how it works in the lowest dimensions, namely for curves and surfaces.

The case of curves is the simplest, as for curves the Kodaira dimension is directly related to positivity properties of the canonical line bundle, hence to the genus $g$ of the curve. Let $X$ be a complex curve, then we have three possibilities: either $X \simeq \IP^1$ in which case $g(X) = 0$, none of the pluricanonical bundles $\omega_X^{\tensor m}$ has holomorphic sections and $\kod(X) = -\infty$, or $X$ is an elliptic curve in which case $g(X) = 1$, all the pluricanonical bundles are trivial and $\kod(X) = 0$, or finally $X$ is a curve of general type, the canonical bundle is ample (hence all the $\omega_X^{\tensor m}$'s do admit holomorphic sections) and $\kod(X) = 1$.

The case of surfaces is more articulated. Since bimeromorphic manifolds have the same Kodaira dimension, we can assume without loss of generality that $X$ is minimal. Smooth minimal complex surfaces were classified by Kodaira (see \cite[Theorem 21 and 22]{Kodaira-I}). Surfaces of general type, i.e., whose canonical bundle is ample are of maximal Kodaira dimension. On the other side, rational surfaces have no non-trivial pluricanonical sections and their Kodaira dimension is $-\infty$. Anyway it is no longer true that surfaces with Kodaira dimension $0$ have trivial canonical bundle: it can be a torsion element in the Picard group of the surface.

The behavior of Kodaira dimension under deformations of the complex structures was also object of study since its introduction. It is known that Kodaira dimension is constant under holomorphic deformations for smooth algebraic manifolds, but this is no longer true if the manifold being deformed is just assumed to be complex.

In the last years all the questions concerning Kodaira dimension have attracted a new attention, as it was generalized by H.~Chen and W.~Zhang in the case of almost complex manifolds (see \cite{ChenZhang-I} and \cite{ChenZhang-II}).

In this paper we will recall their definition and discuss some examples highlighting some features of Kodaira dimension on genuine almost complex manifolds and their deformations. Finally, we make some speculations on the invariance of Kodaira dimension for bimeromorphic almost complex manifolds.

\section{Preliminaries}

All the manifolds we consider throughout this paper are smooth compact differentiable manifolds. Let $M$ be a smooth $2n$-dimensional manifold, as an extra datum we will consider an \emph{almost complex structure} $J$ on $M$, i.e., an endomorphism of the tangent bundle $TM$ of $M$ which squares to $-\id_M$:
\[J: TM \longrightarrow TM, \qquad J^2 = -\id_M.\]

As a consequence $(TM, J)$ becomes a complex bundle on $M$ of rank $n$, which is naturally isomorphic to a subbundle of the complexified tangent bundle $T_\IC M = TM \tensor_\IR \IC$:
\[\begin{array}{ccl}
(TM, J) & \longrightarrow & T_\IC M\\
v \in (TM)_x & \longmapsto & \frac{1}{2}(v - \ii Jv) \in (T_\IC M)_x = (TM)_x \otimes_\IR \IC.
\end{array}\]
More precisely, this map identifies $(TM, J)$ with the eigenbundle of $T_\IC M$ relative of the eigenvalue $\ii$ of the complex extension $J_\IC$ of $J$. We call such eigenbundle $T^{1, 0}_J M$ or $T^{1, 0} X$.

\begin{defin}
An \emph{almost complex manifold} $X$ is a pair $(M, J)$ where $M$ is a $2n$-dimensional smooth differentiable manifold and $J$ is an almost complex structure on $M$.
\end{defin}

The complexified tangent bundle of an almost complex manifold $X = (M, J)$ splits as the direct sum of $T^{1, 0} X$ and its conjugate:
\[T_\IC M = T^{1, 0}_J M \oplus T^{0, 1}_J M, \qquad T^{0, 1}J M = \overline{T^{1, 0}_J M}.\]
This decomposition induces a decomposition both on the complexified cotangent bundle $T^*_\IC M$ and on its exterior powers, so any complex valued smooth $k$-form on $X$ can be written as a sum of forms of type $(p, q)$ with $p + q = k$:
\[\bigwedge^k T^*_\IC X = \bigoplus_{p + q = k} {T^*}^{p, q} X, \qquad {T^*}^{p, q} X = \bigwedge^p {T^*}^{1, 0} X \otimes \bigwedge^q {T^*}^{0, 1} X.\]
Accordingly the exterior differential $d$ splits in the sum of four differential operators:
\[d = \mu + \del + \delbar + \mubar, \qquad \left\{ \begin{array}{l}
\mu: {T^*}^{p, q} X \longrightarrow {T^*}^{p + 2, q - 1} X\\
\del: {T^*}^{p, q} X \longrightarrow {T^*}^{p + 1, q} X\\
\delbar: {T^*}^{p, q} X \longrightarrow {T^*}^{p, q + 1} X\\
\mubar: {T^*}^{p, q} X \longrightarrow {T^*}^{p - 1, q + 2} X.
\end{array} \right.\]

In terms of these operators, it is known that the almost complex manifold $X$ is actually a \emph{complex manifold} if and only if $\mubar \equiv 0$.

\section{Kodaira dimension}

Let $X = (M, J)$ be an almost complex manifold of (complex) dimension $n$. The bundle ${T^*}^{n, 0} X$ is then a complex bundle on $X$ of rank $1$, which is called the \emph{canonical bundle} of $X$ and denoted by $\omega_X$. Its $m^{\text{th}}$ tensor power $\omega_X^{\tensor m}$ is called a \emph{pluricanonical} or \emph{$m$-canonical} bundle.

Let us recall how we can extend the differential operator $\delbar$ to a differential operator $\delbar: \omega_X^{\tensor m} \longrightarrow {T^*}^{0, 1} X \otimes \omega_X^{\otimes m}$. We begin with the case $m = 1$, and take a trivialization for ${T^*}^{n, 0}$ so that any local smooth section $\varphi$ can be written as $\varphi = f \cdot \varphi^1 \wedge \ldots \wedge \varphi^n$ for a suitable smooth complex valued function $f$ and $(1, 0)$-forms $\varphi^1, \ldots, \varphi^n$. Then we define $\delbar \varphi$ applying the usual Leibnitz rule to this expression. This is the extension to the pluricanonical bundles: given a smooth local $m$-canonical form $\phi_1 \tensor \ldots \tensor \phi_m$ we define
\[\delbar\left( \phi_1 \tensor \ldots \tensor \phi_m \right) = \sum_{i = 1}^m \phi_1 \tensor \ldots \tensor \underbrace{\delbar \phi_i}_{i^{\text{th}} \text{ place}} \tensor \ldots \tensor \phi_m .\]

\begin{defin}
Let $X = (M, J)$ be an almost complex manifold. We say that a smooth pluricanonical form $\varphi \in \Gamma(X, \omega_X^{\tensor m})$ is \emph{pseudoholomorphic} if $\delbar \varphi = 0$, and denote
\[H^0(X, \omega_X^{\otimes m}) = \set{\varphi \in \Gamma(X, \omega_X^{\tensor m}) \st \delbar \varphi = 0}.\]
\end{defin}

All what we recalled up to now can be extended to a more general setting, where we consider a smooth complex bundle $E$ on $X$ endowed with a compatible almost complex structure $\cJ$ (a so called \emph{bundle almost complex structure}, see \cite{deBartolomeisTian}). In this case the bundle almost complex structure $\cJ$ uniquely defines an operator $\delbar_E: \Gamma(X, E) \longrightarrow \Gamma\pa{X, {T^*}^{0, 1} X \tensor E}$ and the space $H^0(X, E) = \set{\varphi \in \Gamma(X, E) \st \delbar_E \varphi = 0}$ can be identified with a suitable space of harmonic forms with values in $E$.

\begin{rem}
In general, to define harmonic forms one actually needs a Hermitian metric $h$ on $X$ with respect to which it is possible to define the formal adjoint $\delbar_E^*$ of $\delbar_E$. In our case we can avoid the reference to any metric since a smooth section of $E$ is a $(0, 0)$-form with values in $E$ and applying $\delbar_E^*$ to any $(0, 0)$-form gives $0$ as result by bidegree reasons.
\end{rem}

The possibility to link these spaces with some space of harmonic forms gives us the opportunity to use results from the theory of elliptic operators. In particular, the following theorem is crucial for the definition of Kodaira dimension.

\begin{thm}[{cf.\ \cite[Theorem 3.6]{ChenZhang-I}}]\label{thm: usc plurigenera}
Let $E$ be a complex vector bundle on a smooth Hermitian almost complex manifold $X$, and let $\delbar_E$ be a pseudodifferential operator. Then the spaces $\cH^{p, q}_{\delbar_E}(X, E)$ of harmonic $(p, q)$-forms on $X$ with values in $E$ are finite dimensional for all $(p, q)$.
\end{thm}

In particular this holds true for the spaces $H^0(X, \omega_X^{\tensor m})$, so the following definition is well-posed.

\begin{defin}[{\cite[Definition 4.3]{ChenZhang-I}}]
Let $X$ be a smooth almost complex manifold. The \emph{$m^{\text{th}}$ plurigenus} of $X$ is
\[P_m(X) = \dim H^0(X, \omega_X^{\tensor m}).\]
The \emph{Kodaira dimension} of $X$ is then
\[\kod(X) = \kod(M, J) = \left\{ \begin{array}{ll}
-\infty & \text{if } P_m(X) = 0 \text{ for all } m \geq 1,\\
\limsup_{m \rightarrow +\infty} \frac{\log P_m(X)}{\log m} & \text{otherwise}.
\end{array}\right.\]
\end{defin}

\subsection{Differences with the classic case and methods of computation}

In the complex case there are many tools for computing the Kodaira dimension of a given manifold. We show here some example.

\begin{exam}
Let $X$ be a complex manifold with trivial canonical bundle $\omega_X \simeq \cO_X$. Then all the pluricanonical bundles are trivial as well and so
\[\kod(X) = 0.\]
\end{exam}

\begin{rem}
A very large class of complex manifolds with trivial canonical bundle is provided by Calabi--Yau manifolds. A Calabi--Yau manifold $X$ is a compact K\"ahler manifold with trivial canonical bundle $\cO_X \simeq \omega_X$ and such that $h^{p, 0}(X) = 0$ for $0 < p < \dim_\IC X$. In \cite{deBarTom} de Bartolomeis and Tomassini suggested the following extension (see \cite[Definition 3.1]{deBarTom}): a \emph{generalized Calabi--Yau manifold} is $X = (M, J, \sigma, \varepsilon)$ where $M$ is a compact $2n$-dimensional differentiable manifold, $J$ is an almost complex structure on $M$, $\sigma$ is a symplectic $2$-form and $\varepsilon$ is a smooth canonical section of $(M, J)$, such that
\begin{enumerate}
\item $g_J(\blank, \blank) = \sigma(J \blank, \blank)$ is a positive definite $J$-Hermitian metric on $M$, with Levi-Civita connection $\nabla^{\text{LC}}$;
\item $\displaystyle \varepsilon \wedge \bar{\varepsilon} = (-1)^{\frac{n(n + 1)}{2}} \pa{\ii}^n \frac{\sigma^n}{n!}$;
\item $\nabla^J \varepsilon = 0$, where $\nabla^J = \nabla^{\text{LC}} - \frac{1}{2} J \nabla^{\text{LC}} J$.
\end{enumerate}
If a generalized Calabi--Yau manifold $X$ satisfies $c_1(X) = 0$, then its canonical bundle $\omega_X$ is \emph{differentiably} trivial. The class of generalized Calabi--Yau manifolds whose first Chern class is trivial is then an interesting and natural place where to study and test further properties of Kodaira dimension.
\end{rem}

\begin{exam}
Let $X$ be a complex manifold such that $\omega_X^{-1}$ is ample. Then all the pluricanonical bundles $\omega_X^{\tensor m}$ are anti-ample and so by Kodaira Vanishing Theorem we have that $H^0(X, \omega_X^{\tensor m}) = 0$ for all $m \geq 1$. It follows that $\kod(X) = -\infty$.
\end{exam}

\begin{exam}
Let $X$ be a complex manifold of general type, i.e., whose canonical bundle is ample. In this case we can use Kodaira vanishing theorem to deduce that the $m^{\text{th}}$ plurigenus of $X$ coincides with the holomorphic Euler characteristic of $\omega_X^{\tensor m}$: $P_m(X) = \dim H^0(X, \omega_X^{\tensor m}) = \chi(X, \omega_X^{\tensor m})$. By Hirzebruch--Riemann--Roch this Euler characteristic can be computed, so we have that
\[P_m(X) = \int_X \operatorname{ch}\pa{\omega_X^{\tensor m}} \operatorname{td}(X) = \alpha m^{\dim(X)} + \ldots.\]
As a consequence
\[\kod(X) = \dim(X).\]
\end{exam}

Up to now, if $X$ is an almost complex manifold but not necessarily complex, theorems like Kodaira Vanishing or Hirzebruch--Riemann--Roch are not known to hold, hence the only way to computate the Kodaira dimension is to determine explicitly the spaces $H^0(X, \omega_X^{\tensor m})$.

The methods employed are quite different from the `standard' one used in complex and algebraic geometry, for example in \cite{ChenZhang-I}, \cite{CattaneoNanniciniTomassini-I} and \cite{CattaneoNanniciniTomassini-II} the Kodaira dimension of some families of almost complex structures on $2$- and $3$-dimensional solvmanifolds were computed. Apart from a common starting point, the methods used there can be packed into two classes.
\begin{enumerate}
\item The common starting point is to translate the problem of determining $H^0(X, \omega_X^{\tensor m})$ into a suitable system of differential equations. In the aforementioned papers, smooth sections of $\omega_X^{\tensor m}$ are of the form $f \cdot \psi^{\tensor m}$ where $f$ is a smooth complex valued function and $\psi$ is a nowhere vanishing smooth $(n, 0)$-form. Write $f = u + \ii v$, then the pseudoholomorphic condition $\delbar\pa{f \cdot \psi^{\tensor m}}$ is a differential system for the unknown $u$ and $v$.
\item A first method is to apply results from \emph{elliptic operators of the second order}. Some of the equations are in fact of the type $\bar{V}(f) = 0$ for some vector field $V$ of type $(1, 0)$. This implies that $V\bar{V}(f) = 0$, and $V\bar{V}$ is a real second order operator, which in some cases is elliptic. Hence we get euqations for $u$ and $v$, and in the cases where the operator is elliptic we can deduce that $u$ and $v$ are constant by the Maximum Principle and the fact that the underlying manifold $M$ is compact.
\item A second method is to use \emph{Fourier analysis}. It may happen that the equations arising from the first step are of the form $\bar{V}(f) + g \cdot f = 0$, where $g$ is a smooth complex valued function. In this case the strategy used in the previous point may fail as the possibility of applying the Maximum Principle strongly depends on $g$. As in all the examples discussed the underlying real manifold is a quotient of $\IR^4$ or $\IR^6$ by the action of a suitable lattice $\Gamma$, we can try to solve the system on $\IR^4$, respectively on $\IR^6$ and then find among all the solutions those which are $\Gamma$-periodic. The fact that all the lattices considered in the examples contain a sublattice acting on $\IR^4$ and $\IR^6$ by translations make it possible to expand the solutions in Fourier series. Hence the differential equations to be solved become equations for the Fourier coefficients, which are easier to solve.
\end{enumerate}

\begin{exam}
We outline here how it is possible to use Fourier analysis to compute the Kodaira dimension with an explicit example. Let $G = (\IR^4, *)$ be the Lie group whose multiplication is
\[(a, b, c, d) * (x, y, z, t) = \pa{x + a, y + b, z + ay + c, t + \frac{1}{2}a^2 y + az + d},\]
and consider the subgroup
\[\Gamma = \set{(a, b, c, d) \in G \st a, b, c, d \in 2\IZ}.\]
Define $\cN = \Gamma \backslash G$, then $\cN$ is a compact four dimensional nilmanifold which admits no integrable almost complex structure. Observe that the subgroup
\[\Gamma' = \set{(a, b, c, d) \in \Gamma \st a = 0}\]
is isomorphic to $(2\IZ)^3$ and acts on $\IR^4$ by translations. As a consequence, a smooth function $f: \cN \longrightarrow \IC$ can be seen as a smooth function $f: \IR^4 \longrightarrow \IC$ such that
\[f(x, y, z, t) = f\pa{(2\alpha, 2\beta, 2\gamma, 2\delta)*(x, y, z, t)} \qquad \forall \alpha, \beta, \gamma, \delta \in \IZ.\]
In particular $f(x, y, z, t) = f(x, y + 2\beta, z + 2\gamma, t + 2\delta)$ for all $\beta, \gamma, \delta \in \IZ$, i.e., $f$ is periodic in $y, z, t$ of period $2$ and so we can write (see \cite[Remark 3.4]{CattaneoNanniciniTomassini-II})
\[f(x, y, z, t) = \sum_{I = (a, b, c) \in \IZ^3} f_I(x) e^{\ii \pi (ay + bz + ct)}.\]
The following $1$-forms on $\IR^4$ are $\Gamma$-invariant hence descend to $1$-forms on $\cN$:
\[e^1 = dx, \qquad e^2 = dy, \qquad e^3 = dz - x dy, \qquad e^4 = dt + \frac{1}{2}x^2 dy - x dz,\]
dually we have the following vector fields, which define a parallelism on $\cN$:
\[e_1 = \frac{\partial}{\partial x}, \qquad e_2 = \frac{\partial}{\partial y} + x\frac{\partial}{\partial z} + \frac{1}{2}x^2\frac{\partial}{\partial t}, \qquad e_3 = \frac{\partial}{\partial z} + x\frac{\partial}{\partial t}, \qquad e_4 = \frac{\partial}{\partial t}.\]
Let $J$ be the almost complex structure on $\cN$ defined by
\[J e_1 = e_3, \qquad J e_2 = e_4, \qquad J e_3 = -e_1, \qquad J e_4 = -e_2.\]
From the complex point of view we have then the $(1, 0)$-forms
\[\varphi^1 = e^1 + \ii e^3, \qquad \varphi^2 = e^2 + \ii e^4,\]
and we deduce by a direct computation that
\[\begin{array}{l}
\delbar \varphi^1 = -\frac{1}{4}\ii \varphi^1 \wedge \bar{\varphi}^2 + \frac{1}{4}\ii \varphi^2 \wedge \bar{\varphi}^1,\\
\delbar \varphi^2 = \frac{1}{2} \varphi^1 \wedge \bar{\varphi}^1.
\end{array}\]
As a consequence
\[\delbar(\varphi^1 \wedge \varphi^2) = \underbrace{\frac{1}{4}\ii \bar{\varphi}^2}_{\alpha} \wedge \varphi^1 \wedge \varphi^2,\]
so the smooth $m^{\text{th}}$-canonical section $f \cdot (\varphi^1 \wedge \varphi^2)^{\tensor m}$ is pseudoholomorphic if and only if
\[\left\{ \begin{array}{l}
\bar{\cX}_1(f) = 0\\
\bar{\cX}_2(f) + \frac{1}{4}\ii m f = 0,
\end{array} \right.\]
where $\cX_1, \cX_2$ are the vector fields of type $(1, 0)$ dual to $\varphi^1, \varphi^2$. We write $f = u + \ii v$ and explicit $\bar{\cX}_1$ and $\bar{\cX}_2$ in terms of $e_1, e_2, e_3, e_4$: the previous system is then equivalent to
\[\left\{ \begin{array}{l}
e_1(u) - e_3(v) = 0\\
e_1(v) + e_3(u) = 0\\
2e_2(u) - 2e_4(v) - mv = 0\\
2e_2(v) + 2e_4(v) + mu = 0.
\end{array} \right.\]
In terms of the Fourier coefficients $u_I, v_I$ for $u$ and $v$ respectively, the last two equations lead to the system
\[\left\{ \begin{array}{l}
2\ii \pi \pa{a + bx + \frac{1}{2}cx^2} u_I - (2\ii \pi c + m) v_I = 0\\
(2\ii \pi c + m) u_I + 2\ii \pi \pa{a + bx + \frac{1}{2}cx^2} v_I = 0.
\end{array} \right.\]
The determinant of this system is $-4\pi^2 \pa{a + bx + \frac{1}{2}cx^2}^2 + (m + 2\ii \pi c)^2$, whose imaginary part is $4\pi m c$. As a consequence (recall that $m \geq 1$) if $c \neq 0$ then $u_I$ and $v_I$ vanish identically. Assume then that $c = 0$, in this case the determinant is $-4\pi^2(a + bx)^2 + m^2$, which vanishes only if $a + bx = \pm \frac{m}{2\pi}$. For a fixed index $I = (a, b, 0)$ there is at most one value of $x$ for which the determinant vanishes, hence there is an open dense set where the functions $u_I$ and $v_I$ vanish. Since $u_I$ and $v_I$ are continuous we deduce that they must vanish everywhere. As a consequence $f \cdot (\varphi^1 \wedge \varphi^2)^{\tensor m}$ is pseudoholomorphic if and only if $f  0$, which by definition means that
\[\kod(\cN, J) = -\infty.\]
\end{exam}

\begin{rem}
Observe that the almost complex structure $J$ defined in the previous example is similar to the one called $J_1$ in \cite[$\S$7]{CattaneoNanniciniTomassini-II}, but it is not difficult to show that $J$ is not an almost complex structure in the twistor sphere defined there.
\end{rem}

\section{Behavior under deformations}

In this section we want to give a brief account on the properties of Kodaira dimension when the almost complex manifold in consideration is allowed to deform. In this direction the most general result is the following, which follows from the semicontinuity of the dimension of the kernel of elliptic operators as shown in \cite[$\S$4.4]{KodairaMorrow}.

\begin{thm}[{cf.\ \cite[Proposition 6.3]{ChenZhang-I}}]
Let $\Delta = \set{t \in \IC \st |t| > \epsilon}$, and let $X_t = (M, J(t))$ be a smooth family of almost complex manifolds for $t \in \Delta$. Then $P_m(X_t)$ is an upper semicontinuous function of $t$.
\end{thm}

Related questions are the dependence on $t$ of the function $\kod(X_t)$ and the possibility that for some suitable classes of manifolds the plurigenera or at least the Kodaira dimension remain constant in families. In fact, Iitaka introduced the concept of Kodaira dimansion (which ha called \emph{canonical dimension}) in \cite{Iitaka-Ddimension} and he also completely solved this problem in the case of surfaces.

\begin{thm}[{\cite[Theorem III]{Iitaka-Deformationsofsurfaces}}]
Plurigenera of compact complex surfaces are invariant under arbitrary holomorphic deformations.
\end{thm}

As a consequence, the Kodaira dimension is also constant. In the case of higher dimensional manifolds it was later shown by Nakamura (see \cite[Theorem 2]{Nakamura}) that for a general holomorphic deformation of a compact complex manifold, the Kodaira dimension of the fibres of the family, as well as the plurigenera, can be non constant. We recall briefly his construction. Nakamura considered the three dimensional complex solvmanifold $N$ with structure equations
\[d\varphi_1 = 0, \qquad d\varphi_2 = \varphi_1 \wedge \varphi_2, \qquad d\varphi_3 = -\varphi_1 \wedge \varphi_3,\]
where $\set{\varphi_1, \varphi_2, \varphi_3}$ is a global coframe of holomorphic invariant $1$-forms. Call $\set{\vartheta_1, \vartheta_2, \vartheta_3}$ the corresponding dual frame of $(1, 0)$-vector fields, and define
\[\varphi^*_1 = \bar{\varphi}_1, \qquad \varphi^*_2 = e^{z_1 - \bar{z}_1} \bar{\varphi}_2, \qquad \varphi^*_3 = e^{-z_1 + \bar{z}_1} \bar{\varphi}_3.\]
Then he took the deformation of $N$ corresponding to
\[\psi = [(t_1 \vartheta_1 + t_2 \vartheta_2 + t_3 \vartheta_3) \tensor \varphi^*_2 + t_4 \vartheta_3 \tensor \varphi^*_3] \in H^{0, 1}_{\delbar}(N, T^{1, 0} N)\]
where $(t_1, t_2, t_3, t_4) \in \IC^4$ with $||(t_1, t_2, t_3, t_4)|| < \varepsilon$ and observed that
\[P_m(N_t) = \left\{ \begin{array}{ll}
1 & \text{if } (t_1, t_2, t_3, t_4) = (0, 0, 0, 0)\\
0 & \text{if } t_1 \neq 0,
\end{array} \right.\]
from which
\[\kod(N_t) = \left\{ \begin{array}{ll}
0 & \text{if } (t_1, t_2, t_3, t_4) = (0, 0, 0, 0)\\
-\infty & \text{if } t_1 \neq 0.
\end{array} \right.\]

Nevertheless, if one restricts to the realm of complex projective algebraic manifolds it is still true that plurigenera and Kodaira dimension are invariant under deformations, as it was later proven by Siu.

\begin{thm}[{\cite[Corollary 0.3]{Siu}}]
Let $\pi: X \longrightarrow \Delta$ be a holomorphic family of compact complex projective algebraic manifolds over the open $1$-disc $\Delta$ with fibre $X_t$. Let $m$ be any positive integer. Then the complex dimension of $\Gamma(X_t, \omega_{X_t}^{\tensor m})$ is independent of $t \in \Delta$.
\end{thm}

In the final part of this section we focus on the behavior of Kodaira dimension for deformations of almost complex manifolds, and show that this is not invariant for generic pseudoholomorphic families.

\begin{rem}
The fact that plurigenera and Kodaira dimension are not invariant under deformations of the almost complex structure was already pointed out with many examples in \cite{ChenZhang-I}, \cite{CattaneoNanniciniTomassini-I} and \cite{CattaneoNanniciniTomassini-II}. All the examples in these papers consider \emph{smooth} deformations of the almost complex structures, depending smoothly on some \emph{real} parameters. The example we present now shows that Kodaira dimension can have a very odd behavior even if we consider a \emph{pseudoholomorphic} deformation of the almost complex structure.
\end{rem}

\begin{exam}
Consider the group of matrices
\[G = \set{\left( \begin{array}{ccccc}
1 & x & z & 0 & 0\\
0 & 1 & y & 0 & 0\\
0 & 0 & 1 & 0 & 0\\
0 & 0 & 0 & 1 & t\\
0 & 0 & 0 & 0 & 1
\end{array}\right) \st x, y, z, t \in \IR},\]
and the action by left multiplication by elements of the lattice $\Gamma$ given by the matrices in $G$ with integral entries. The quotient $\Gamma \backslash G$ is then a smooth differentiable manifold $M$, wich is the underlying real manifold of the Kodaira--Thurston surface (see \cite{Thurston} or \cite[$\S$6]{Kodaira-I}).
Let
\[e_1 = \partial_t, \qquad e_2 = \partial_x, \qquad e_3 = \partial_y + x \partial_z, \qquad e_4 = \partial_z:\]
this is a global frame of $\Gamma$-invariant vector fields on $G$ which descends to a global frame of vector fields on $M$ (we will keep the same notation for the induced fields on $M$). Let now $\Delta = \set{t \in \IC \st |t| < \pi}$ and
\[e_5 = \partial_t + \partial_{\bar{t}}, \qquad e_6 = \ii \pa{\partial_t - \partial_{\bar{t}}},\]
be the usual derivatives with respect to the real and imaginary part of $t$. Define the almost complex structure on $M \times \Delta$ described in the global frame $\set{e_1, e_2, e_3, e_4, e_5, e_6}$ by the matrix
\[J = \left( \begin{array}{cccccc}
0 & -1 & 0 & 0 & 0 & 0\\
1 & 0 & 0 & 0 & 0 & 0\\
0 & 0 & -\frac{\Im(t)}{\Re(t) + \pi} & -\frac{|t + \pi|^2}{\Re(t) + \pi} & 0 & 0\\
0 & 0 & \frac{1}{\Re(t) + \pi} & \frac{\Im(t)}{\Re(t) + \pi} & 0 & 0\\
0 & 0 & 0 & 0 & 0 & -1\\
0 & 0 & 0 & 0 & 1 & 0
\end{array} \right).\]
With this almost complex structure the manifold $M \times \Delta$ is an almost complex manifold, as well as every manifold $M_t = M \times \set{t}$ which is an almost complex manifold by means of the almost complex structure $J_t$. Consider now the projection $\pi: M \times \Delta \longrightarrow \Delta$. With respect to the frames $\set{e_1, e_2, e_3, e_4, e_5, e_6}$ and $\set{e_5, e_6}$ on $M \times \Delta$ and $\Delta$ respectively the differential $d\pi$ is described by the matrix
\[d\pi = \left( \begin{array}{cccccc}
0 & 0 & 0 & 0 & 1 & 0\\
0 & 0 & 0 & 0 & 0 & 1
\end{array} \right)\]
and it is then easy to see that $\pi$ is pseudoholomorphic (of course, we consider on $\Delta$ the standard complex structure). Assume now that $t$ is \emph{real}: in this case the almost complex structure on $M_t$ is
\[J_t = \left( \begin{array}{cccc}
0 & -1 & 0 & 0\\
1 & 0 & 0 & 0\\
0 & 0 & 0 & -(t + \pi)\\
0 & 0 & \frac{1}{t + \pi} & 0
\end{array} \right),\]
which is, up to a translation, the almost complex structure described in \cite[$\S$6.1]{ChenZhang-I}. By \cite[Proposition 6.1]{ChenZhang-I} we can then conclude that
\[\kod(M_t, J_t) = \left\{ \begin{array}{ll}
-\infty & \text{if } t \in \Delta \cap \pa{\IR \smallsetminus \pi \IQ},\\
0 & \text{if } t \in \Delta \cap \pi \IQ.
\end{array} \right.\]
Observe that while it is still true by Theorem \ref{thm: usc plurigenera} that each plurigenus is independently upper semicontinuous as a function of  $t$, on the contrary there is no semicontinuity property for the Kodaira dimension (not even if one considers only pseudoholomorphic deformations)
\end{exam}

\section{A speculation on meromorphic functions on almost complex manifolds}

In this final section we want to discuss about the invariance of Kodaira dimansion for bimeromorphic manifolds and the links with the canonical divisor. In particular, we want to give some directions for further studies on the geometry of almost complex manifolds from a bimeromorphic point of view.

Let us consider the complex case for a moment. The canonical bundle $\omega_X$ of a complex manifold $X$ is a holomorphic line bundle, and if we let $f$ be a meromorphic section of $\omega_X$ then it is possible to define a \emph{canonical divisor} on $X$ as the divisor $(f)$ associated to $f$. Moreover, if $X$ and $X'$ are bimeromorphic complex manifolds then $\kod(X) = \kod(X')$. This shows that the concept of \emph{meromorphic functions} or \emph{sections} is central in complex geometry.

In the almost complex setting it is not clear how to define meromorphic functions and sections. Here we give some ideas on how these concepts may be defined.

On an almost complex manifold we can naturally define a sheaf of rings, which associates to an open subset $U \subseteq X$ the ring
\[\cO_X(U) = \set{f: U \longrightarrow \IC \st f \text{ is pseudoholomorphic}}.\]
Hence $(X, \cO_X)$ is a locally ringed space. As such we have the possibility of defining meromorphic functions in a quite straightforward way, as explained in \cite[$\S$20.1]{EGA-IV-IV}.

Let $U \subseteq X$ be an open subset, and let $S(U)$ be the subset of $\cO_X(U)$ consisting of the holomorphic functions which are not zero-divisors. Then we define
\[\cM_X(U) = S(U)^{-1} \cO_X(U),\]
the localization of $\cO_X(U)$ with respect to the multiplicative system $S(U)$. A meromorphic function on $U$ is then an element of $\cM_X(U)$.

It is possible to show that $\cM_X(U)$ is a local ring, with maximal ideal given by
\[\set{\frac{f}{g} \st f \text{ is a zero-divisor}, g \text{ is not a zero-divisor}}\]
and that the canonical map
\[\cO_X(U) \longrightarrow \cM_X(U)\]
is injective.

It becomes then of interest to determine which pseudoholomorphic functions are zero-divisors. With respect to this problem, we observe that almost complex manifolds sit between two different extremes. In the case of smooth differentiable manifolds a smooth function $f: U \longrightarrow \IR$ is a zero-divisor if and only if $f^{-1}(0)$ ha non-empty interior. On the contrary on a complex manifold $X$ we have that $\cO_X(U)$ is an integral domain, hence the only zero-divisor is the zero function. The study of the zero locus of a pseudoholomorphic function is then crucial to understand the ring $\cM_X(U)$.

The following step is to consider meromorphic sections of line bundles, i.e., those sections which are given by meromorphic functions on local trivializations of the line bundle. In this case the knowledge of the structure of the zero locus of pseudoholomorphic functions is needed to understand if it is possible to define a divisor from a meromorphic section.

Finally, we make a remark on meromorphic maps $X \longrightarrow X'$. On a complex manifold it is possible to define holomorphic charts, which exhibits it locally as an open set of $\IC^n$. Hence it is possible to define a meromorphic map from complex manifolds as a map which is meromorphic when we write its expression in coordinates. On the contrary, on an almost complex manifold we do not have the possibility of speaking of complex charts, but all we have are the real charts and the almost complex structure. This shows the main difficulty to give a definition of meromorphic map between almost complex manifolds.

Of course, a good definition of meromorphic function between almost complex manifolds should enjoy the property that if $X$ and $X'$ are bimeronorphic almost complex manifolds then $\kod(X) = \kod(X')$.

\bibliographystyle{plain}
\bibliography{Biblio}

\end{document}